# RATNER'S THEOREM AND INVARIANT THEORY

Akihiko Yukie

ABSTRACT. In this note, we consider applications of Ratner's theorem to constructions of families of polynomials with dense values on the set of primitive integer points from the viewpoint of invariant theory.

Let $Q$ be a non-degenerate indefinite quadratic form in $n \geq 3$ variables. We say $Q$ is irrational if the ratio of at least one pair of coefficients is irrational.

The following theorem was called the Oppenheim conjecture and was proved in final form by Margulis [9].

**Theorem 1 (Margulis)** *If $Q(x)$ is an irrational non-degenerate indefinite quadratic form in $n \geq 3$ variables, the set $\{Q(x) \mid x \in \mathbb{Z}^n\}$ is dense in $\mathbb{R}$.*

Many partial results were known prior to the work of Margulis (for example, Oppenheim [10], [11], [12], Davenport (with others) [6], [4] [5], [1], [19]). Note that because of a result of Lewis [7], these results imply the Oppenheim conjecture for such $Q$ with $n \geq 21$. For forms of higher degree, there is a result by Schmidt [22].

Raghunathan reduced the Oppenheim conjecture to the following problem.

**Problem 2** Show that $\mathrm{SO}(Q)_{\mathbb{R}} \mathrm{SL}(3)_{\mathbb{Z}}$ is dense in $\mathrm{SL}(3)_{\mathbb{R}}$.

What Margulis did was to solve the above problem, and in fact he proved that the values of $Q$ on the set of primitive integer points are dense in $\mathbb{R}$.

This problem can be viewed as a problem on orbit closures in Lie groups. We first recall Ratner's theorem.

**Theorem 3 (Ratner)** *Let $G$ be a connected Lie group and $U$ a connected subgroup of $G$ generated by unipotent elements of $G$. Then given any lattice $\Gamma \subset G$ and $x \in G$, there exists a connected closed subgroup $U \subset F \subset G$ such that $\overline{Ux\Gamma} = Fx\Gamma$. Moreover, $F/(F \cap (x\Gamma x^{-1}))$ has a finite invariant measure.*

The first statement was called Raghunathan's topological conjecture, and the second statement was proved by Ratner in conjunction with Raghunathan's topological conjecture. Raghunathan's topological conjecture was published by Dani [3] for one dimensional unipotent groups and was generalized to groups generated by unipotent elements by Margulis [8]. The proof for the general case was given by

1991 *Mathematics Subject Classification*. Primary 11J25; Secondary 22E40.
Partially supported by NSF grant DMS-9401391

Ratner in a series of papers [14], [15], [16], [17]. For these, there is an excellent survey article by Ratner [18].

Note that $F$ in the above theorem is unique. To see this, one should note that the image of $F$ in $G/\Gamma$ by the map $g \to gx\Gamma$ is $F/(F \cap (x\Gamma x^{-1}))$ and $F \cap (x\Gamma x^{-1})$ is a discrete subgroup of $F$. So the tangent space of $F$ at the unit element is determined by the image of $F$ in $G/\Gamma$. Also if $\Gamma_1, \Gamma_2$ are commensurable lattices and $\overline{Ux\Gamma_i} = F_i x\Gamma$ for $i = 1, 2$, $F_1 = F_2$. To see this, one can assume that $\Gamma_1 \subset \Gamma_2$ and the index of $\Gamma_1$ in $\Gamma_2$ is finite. Then if $g_1, \cdots, g_l$ are representatives of $\Gamma_1 \backslash \Gamma_2$,

$$\overline{Ux\Gamma_2} = \cup_{i=1}^l \overline{Ux\Gamma_1} g_i = \cup_{i=1}^l F_1 x\Gamma_1 g_i = F_1 x\Gamma_2.$$

Because of the uniqueness of $F$, $F_2 = F_1$.

It is easy to deduce Problem 2 (and hence Theorem 1) from Ratner's theorem. We will not give the proof here, but we will use a similar argument later in Theorems 7,14 to construct a family of quartic forms in four variables whose values on the set of primitive integer points are dense in $\mathbb{R}$.

For any algebraic group $G$ defined over a subfield of $\mathbb{R}$, let $G^\circ_{\mathbb{R}+}$ be the identity component of $G_\mathbb{R}$ in the classical topology. We denote the identity component of $G$ in the Zariski topology by $G^\circ$. Therefore, $G^\circ_\mathbb{R}$ is the set of $\mathbb{R}$-rational points in $G^\circ$, and this group may not be connected in the classical topology. If $G$ is an algebraic group over $\mathbb{C}$, it is connected in the Zariski topology if and only if $G_\mathbb{C}$ is connected in the classical topology. So the identity component of $G_\mathbb{C}$ in the classical topology is an algebraic group over $\mathbb{C}$. We use the notation $G^\circ_\mathbb{C}$ for this identity component and this can be regarded as the set of $\mathbb{C}$-valued points in the identity component $G^\circ$ in the Zariski topology. If $F$ is an algebraic group over $\mathbb{C}$ and $k$ is a subfield of $\mathbb{C}$, we say $F$ is defined over $k$ if there exists an algebraic group $G$ over $k$ such that $F = G_\mathbb{C}$. Suppose $G$ is an algebraic group defined over a subfield $k$ of $\mathbb{C}$. Then it is known that the identity component $G^\circ_\mathbb{C}$ is defined over $k$ (see [2, p. 47]).

If a connected Lie group $F$ is in the form $F = G^\circ_{\mathbb{R}+}$ for an algebraic group $G$ over $\mathbb{Q}$, we say that $F$ is defined over $\mathbb{Q}$. Moreover, we use the notation $F_\mathbb{C}$ for $G_\mathbb{C}$ even though this is a slight abuse of notation. By the above remark, $G^\circ_\mathbb{C}$ is defined over $\mathbb{Q}$. So we may choose such $G$ so that $G_\mathbb{C}$ is connected. This $G$ is automatically connected in the Zariski topology over any extension of $\mathbb{Q}$.

Now let $H_1 \subset H_2$ be algebraic groups over $\mathbb{Q}$ (and so $H_1$ is a closed subgroup). We assume that $H_2$ is semi-simple, and $H_{1\mathbb{C}}, H_{2\mathbb{C}}$ are connected. Consider $H^\circ_{1\mathbb{R}+}, H^\circ_{2\mathbb{R}+}$. For any closed connected Lie subgroup $H \subset H^\circ_{2\mathbb{R}+}$, we consider the following condition.

**Condition 4** $H$ is generated by unipotent elements. Moreover, the projection of $H$ to any simple component of $H^\circ_{2\mathbb{R}+}$ is non-trivial.

Note that if a connected Lie group $H$ is simple and non-compact, $H$ is generated by unipotent elements. We assume that $H^\circ_{1\mathbb{R}+}$ satisfies Condition 4. Suppose we can carry out the following problem.

**Problem 5** (a) Classify all the closed connected complex Lie subgroups $F$ such that $H_{1\mathbb{C}} \subset F \subset H_{2\mathbb{C}}$ and show that there are only finitely many such $F$'s. Moreover, show that all such $F$'s are defined over $\mathbb{Q}$.

(b) For each $F \neq H_{2\mathbb{C}}$, construct an $H_2$-variety $X_F$ defined over $\mathbb{Q}$ with the property that $F$ has a unique fixed point $w_F \in X_{F\mathbb{Q}}$ in $X_{F\mathbb{C}}$.

Of course, we may not always be able to carry out Problem 5. Note that the existence of $X_F$ in (2) in many cases follows from a theorem of Borel and Harish–Chandra. More precisely, suppose $F$ in Problem 5 satisfies the following condition

$$\{h \in H_{2\mathbb{C}} \mid hfh^{-1} \in \mathrm{N}_{H_{2\mathbb{C}}}(F) \text{ for all } f \in F\} = \mathrm{N}_{H_{2\mathbb{C}}}(F)$$

($\mathrm{N}_{H_{2\mathbb{C}}}(F)$ is the normalizer of $F$ in $H_{2\mathbb{C}}$). Then there exists a faithful representation $H_2 \to \mathrm{GL}(E)$ of $H_2$ defined over $\mathbb{Q}$ and a point $v \in \mathbb{P}(E)_{\mathbb{Q}}$ such that $\mathrm{N}_{H_{2\mathbb{C}}}(F)$ is the stabilizer of $v$ in $H_{2\mathbb{C}}$ (see 5.1 Theorem [2, p. 89]). Let $X$ be the image of $H_2$ in $\mathbb{P}(E)$ by the map $H_2 \ni h \to hv \in \mathbb{P}(E)$. Then $X$ is a locally closed subset of $\mathbb{P}(E)$ defined over $\mathbb{Q}$ (since $X$ is a homogeneous space). Then $y = hv \in X_{\mathbb{C}}$ ($h \in H_{2\mathbb{C}}$) is a fixed point of $F$ if and only if $hfh^{-1} \in \mathrm{N}_{H_{2\mathbb{C}}}(F)$ for all $f \in F$, and so $h \in \mathrm{N}_{H_{2\mathbb{C}}}(F)$. So we may choose this $X$ as $X_F$. We will consider Problem 5 for certain situations (see Problem 9), and so far the above condition is satisfied for all the cases we considered.

In the following, we assume that we carried out Problem 5.

**Definition 6** An element $g \in H_{2\mathbb{C}}$ is *sufficiently irrational* if $gw_F \notin X_{F\mathbb{Q}}$ for all $F \neq H_{2\mathbb{C}}$.

Note that if $F_1, F_2$ are such groups and there exists an $H_2$-equivariant map $X_{F_1} \to X_{F_2}$ defined over $\mathbb{Q}$, we don't have to consider $X_{F_1}$, because $X_{F_1\mathbb{Q}}$ maps to $X_{F_2\mathbb{Q}}$. Let $\Gamma \subset H_{2\mathbb{R}+}^\circ$ be an arithmetic lattice. Now the following theorem is an easy consequence of Ratner's theorem.

**Theorem 7** *Suppose $g \in H_{2\mathbb{C}}$ is sufficiently irrational, $gH_{1\mathbb{C}}g^{-1}$ is defined over $\mathbb{R}$, and the identity component $U$ of $(gH_{1\mathbb{C}}g^{-1}) \cap H_{2\mathbb{R}}$ satisfies* Condition 4. *Then $U\Gamma$ is dense in $H_{2\mathbb{R}+}^\circ$.*

*Proof.* By Ratner's theorem, there exists a closed connected subgroup $U \subset F_1 \subset H_{2\mathbb{R}+}^\circ$ such that $\overline{U\Gamma} = F_1\Gamma$. Moreover, $F_1$ is defined over $\mathbb{Q}$ by Proposition (3.2) [23, pp. 321–322]. Note that Shah considered a particular lattice in [23], but by our commensurability discussion of lattices after Ratner's theorem, we can consider any arithmetic lattice.

Suppose $F_1 \neq H_{2\mathbb{R}+}^\circ$. Then $\dim F_1 < \dim H_{2\mathbb{R}+}^\circ$. Since both groups are defined over $\mathbb{Q}$, $F_{1\mathbb{C}} \neq H_{2\mathbb{C}}$ by considering the Lie algebras. Since $gH_{1\mathbb{C}}g^{-1}$ is defined over $\mathbb{R}$, $gH_{1\mathbb{C}}g^{-1} \subset F_{1\mathbb{C}}$. So $H_{1\mathbb{C}} \subset g^{-1}F_{1\mathbb{C}}g \neq H_{2\mathbb{C}}$. Since we are assuming we carried out Problem 5, there exists an algebraic group $F_2$ over $\mathbb{Q}$ such that $g^{-1}F_{1\mathbb{C}}g = F_{2\mathbb{C}}$. Since $w_{F_2}$ is the unique fixed point by $g^{-1}F_{1\mathbb{C}}g$, $p_{F_1} = gw_{F_2}$ is the unique fixed point by $F_{1\mathbb{C}}$. Let $\sigma \in \mathrm{Aut}(\mathbb{C}/\mathbb{Q})$. Then $F_{1\mathbb{C}} = F_{1\mathbb{C}}^\sigma$ clearly fixes $p_{F_1}^\sigma$. By the uniqueness of the fixed point, $p_{F_1}^\sigma = p_{F_1}$. Since this is true for all $\sigma \in \mathrm{Aut}(\mathbb{C}/\mathbb{Q})$, $p_{F_1} = gw_{F_2} \in X_{F_2\mathbb{Q}}$, which is a contradiction since we assumed that $g$ is sufficiently irrational. $\square$

We want to use the above observation in the setting of prehomogeneous vector spaces. We now recall the definition of prehomogeneous vector spaces. Let $G$ be a reductive group, $V$ a representation of $G$, and $\chi$ a non-trivial character, all defined over a field $k$ of characteristic zero. We assume that $G_{\bar k}$ is connected in the Zariski topology.

**Definition 8** $(G, V, \chi)$ is called a *prehomogeneous vector space* if
(1) there exists a Zariski open orbit,
(2) there exists a non-zero polynomial $\Delta(x) \in k[V]$ such that $\Delta(gx) = \chi(g)\Delta(x)$ for all $g \in G$, $x \in V$.

Let $V^{\mathrm{ss}} = \{x \in V \mid \Delta(x) \neq 0\}$ and we call it the set of semi-stable points. A polynomial $\Delta(x)$ in (2) is called a relative invariant polynomial. If the representation is irreducible, the choice of $\chi$ is essentially unique and we may write $(G, V)$ also.

For the rest of this paper, we only consider irreducible prehomogeneous vector spaces. Let $(G, V)$ be a prehomogeneous vector space. If there exists a relative invariant polynomial $\Delta(x)$ whose Hessian is not identically zero, $(G, V)$ is called regular. It is known (see [21]) that the regularity is equivalent to the existence of a point $w$ in the open orbit in $V_k^{\mathrm{ss}}$ whose stabilizer $G_w$ is reductive (since ch $k = 0$). An important consequence of the regularity is that $V_{\bar k}^{\mathrm{ss}}$ is a single $G_{\bar k}$-orbit. Up to a certain notion of equivalence, irreducible regular prehomogeneous vector spaces over any algebraically closed field $k$ of characteristic zero were classified by M. Sato and Kimura in [21]. Each equivalence class contains a prehomogeneous vector space $(G, V)$ such that $\dim V$ is the smallest, and such $(G, V)$ is called reduced. Over local and global fields, the classification was carried out by H. Saito [20]. We list the split forms of irreducible reduced regular prehomogeneous vector spaces as follows.

(1) $G = \mathrm{GL}(n) \times G'$, $V = \mathrm{M}(n, n)$ where $G' \subset \mathrm{GL}(n)$ is a reductive subgroup such that $k^n$ is an irreducible representation of $G'$.
(2) $G = \mathrm{GL}(1) \times \mathrm{GL}(n)$, $V = \mathrm{Sym}^2 k^n$.
(3) $G = \mathrm{GL}(2n)$, $V = \wedge^2 k^n$.
(4) $G = \mathrm{GL}(1) \times \mathrm{GL}(2)$, $V = \mathrm{Sym}^2 k^3$.
(5), (6), (7) $G = \mathrm{GL}(1) \times \mathrm{GL}(n)$, $V = \wedge^3 k^n$ where $n = 6, 7, 8$ (as in $E_6, E_7, E_8$).
(8) $G = \mathrm{GL}(3) \times \mathrm{GL}(2)$, $V = \mathrm{Sym}^3 \otimes k^2$.
(9) $G = \mathrm{GL}(6) \times \mathrm{GL}(2)$, $V = \wedge^2 k^6 \otimes k^2$.
(10), (11) $G = \mathrm{GL}(n) \times \mathrm{GL}(5)$, $V = k^n \otimes \wedge^2 k^5$ where $n = 3, 4$.
(12) $G = \mathrm{GL}(3) \times \mathrm{GL}(3) \times \mathrm{GL}(2)$, $V = \mathrm{M}(3, 3) \otimes k^2$.
(13) $G = \mathrm{GSp}(2n) \times \mathrm{GL}(2m)$, $V = k^{2n} \otimes k^{2m}$ where $n \geq 2m \geq 2$.
(14) $G = \mathrm{GSp}(6)$ and $V = k^{14}$ is the irreducible representation which corresponds to the third fundamental weight.
(15) $G = \mathrm{GO}(n) \times \mathrm{GL}(m)$, $V = k^n \otimes k^m$.
(16), (17), (18) $G = \mathrm{GSpin}(7) \times \mathrm{GL}(n)$, $V = k^8 \otimes k^n$ where $k^8$ is the spin representation and $n = 1, 2, 3$.
(19), (22) $G = \mathrm{GSpin}(n)$ where $n = 9, 11$ and $V$ is the spin representation.
(20), (21) $G = \mathrm{GSpin}(10) \times \mathrm{GL}(n)$, $V = k^{16} \otimes k^n$ where $k^{16}$ is the half spin representation and $n = 2, 3$.
(23), (24) $G = \mathrm{GSpin}(n)$ where $n = 12, 14$ and $V$ is the half spin representation.
(25), (26) $G = \mathrm{G}_2 \times \mathrm{GL}(n)$, $V = k^7 \otimes k^n$ where $\mathrm{G}_2$ is the automorphism group of the split octonion, $k^7$ is the imaginary part of the split octonion, and $n = 1, 2$.

(27), (28) $G = \mathrm{E}_6 \times \mathrm{GL}(n)$, $V = J \otimes k^n$ where $\mathrm{E}_6$ is a semi-simple group of type $\mathrm{E}_6$ and $J$ is the exceptional Jordan algebra.
(29) $G$ is a group whose semi-simple part is a group of type $\mathrm{E}_7$ and $V = k \oplus k \oplus J \oplus J$.

In the above list, $\mathrm{GSp}(2n), \mathrm{GO}(n), \mathrm{GSpin}(n)$ are the groups defined similarly as $\mathrm{Sp}(2n), \mathrm{O}(n), \mathrm{Spin}(n)$ except for allowing scalar multiplications.

Now let $(G, V)$ be an irreducible prehomogeneous vector space in (1)–(29). Let $w \in V_{\mathbb{Z}}^{\mathrm{ss}}$, and $H = [G, G]$. Then $H$ is a semi-simple group. Let $H_w$ be the stabilizer of $w$ in $H$. We consider the following two problems.

**Problem 9** (1) Carry out Problem 5 for $H_1 = H_w^\circ$, $H_2 = H$.
(2) Carry out Problem 5 for the image $H_1$ of $H$ in $\mathrm{SL}(V)$ and $H_2 = \mathrm{SL}(V)$.

For Problem 9(2), we do not necessarily have to consider prehomogeneous vector spaces. Let $V$ be a rational representation of a semi-simple algebraic group $H$ over $\mathbb{Q}$ such that $H_{\mathbb{C}}$ is connected. Then we can consider Problem 9(2) for $H_1, H_2$ defined in the same manner.

Consider problem 9(1) for the case (1) in the classification. Let $H' = [G', G']$. Then $H = \mathrm{SL}(n) \times H'$ in this case, and $w = I_n$ (the identity matrix) belongs to $V_{\mathbb{Z}}^{\mathrm{ss}}$. It is easy to see that $H_w = \{(h, {}^t h^{-1}) \mid h \in H'\}$. So Problem 9(1) for this case is almost equivalent to Problem 9(2) for irreducible representations of semi-simple groups. Since this case contains a rather large class of representations, we do not consider the case (1) in the classification when we consider Problem 9(1), and treat it separately as Problem 9(2).

Consider Problem 9(1). If $x \in V_{\mathbb{R}}^{\mathrm{ss}}$, there exists $g \in G_{\mathbb{C}}$ such that $gw = x$. We can choose $t$ in the center of $G_{\mathbb{C}}$ and $h \in H_{\mathbb{C}}$ so that $g = th$.

**Definition 10** A point $x \in V_{\mathbb{R}}^{\mathrm{ss}}$ is *sufficiently irrational* if there exists an element $g = th \in G_{\mathbb{C}}$ as above such that $h$ is sufficiently irrational.

Problem 9(1) is not applicable to cases (2) $n = 2$, (4), (8), (9), (11), (12), (15) $m = 2$, (17), (26), (28) because stabilizers of points in $V^{\mathrm{ss}}$ do not satisfy Condition 4. In [26], [25], [27], we considered cases (3), (5), (6), (7), (10). In [28], we consider all the other applicable cases except for the cases (1), (29). For the case (29), we seem to need a new interpretation of the "Freudenthal quartic". For Problem 9(1), we also construct $G$-equivariant maps $\phi_F : V^{\mathrm{ss}} \to X_F$ for all $F$'s in Problem 5 so that $\phi_F(w) = w_F$. Moreover, the center of $G$ acts trivially on $X_F$ for all $F$. So an alternative way of stating the above definition is that $x \in V_{\mathbb{R}}^{\mathrm{ss}}$ is sufficiently irrational if $\phi_F(x) \notin X_{F\mathbb{Q}}$ for all $F$, because $\phi_F(x) = t\phi_F(w) = tw_F$.

If $x \in V_{\mathbb{R}}^{\mathrm{ss}}$, there exists $g \in G_{\mathbb{C}}$ such that $x = gw$. We write $g = th$ as above. Then $H_{x\mathbb{C}}^\circ$ is defined over $\mathbb{R}$. Since $H_{x\mathbb{C}}^\circ = h H_{w\mathbb{C}}^\circ h^{-1}$, this group is also defined over $\mathbb{R}$. Note that we cannot assume $g \in G_{\mathbb{R}}$ in Definition 10 since $V_{\mathbb{R}}^{\mathrm{ss}}$ is not necessarily a single $G_{\mathbb{R}}$-orbit.

We state an application of Problem 9(1). Consider the case (5) in the classification. Let $e_1, \cdots, e_6$ be the standard basis of $k^6$, and $w = e_1 \wedge e_2 \wedge e_3 + e_4 \wedge e_5 \wedge e_6$. Let $E_1, E_2$ be subspaces of $k^6$ spanned by $\{e_1, e_2, e_3\}, \{e_4, e_5, e_6\}$ respectively. It is known that $w \in V_{\mathbb{R}}^{\mathrm{ss}}$. Let $\mathrm{Gr}(3, 6)$ be the Grassmann of 3-spaces in $k^6$, and $X = (\mathbb{Z}/2\mathbb{Z}) \backslash (\mathrm{Gr}(3, 6) \times \mathrm{Gr}(3, 6))$ where $\mathbb{Z}/2\mathbb{Z}$ acts by exchanging two factors. Let $[E_1], [E_2] \in \mathrm{Gr}(3, 6), ([E_1], [E_2]) \in X$ be the points determined by $E_1, E_2$. Let $x = \sum x_{ijk} e_i \wedge e_j \wedge e_k = gw \in V_{\mathbb{R}}^{\mathrm{ss}}$ where $g \in \mathrm{GL}(6)_{\mathbb{C}}$. Then it turns out that $x$ is sufficiently irrational if $g[E_1], g[E_2] \notin \mathrm{Gr}(3,6)_{\mathbb{Q}}$ and $g([E_1],[E_2]) \notin X_{\mathbb{Q}}$. Both $\mathrm{Gr}(3,6)$ and $X$ are defined over $\mathbb{Q}$ and one can determine whether or not a point in these spaces is rational by considering Plücker coordinates. In this case, $H^\circ_{x\mathbb{R}+}$ satisfies Condition 4 for all $x \in V^{\mathrm{ss}}_{\mathbb{R}}$ [25]. The following result is proved in [25].

**Theorem 11** *Suppose $x \in V^{\mathrm{ss}}_{\mathbb{R}}$ is sufficiently irrational. Then for any $y = (y_{ijk}) \in \wedge^3 \mathbb{R}^5$ and $\epsilon > 0$, there exists a $\mathbb{Z}$-basis $\{u_1, \cdots, u_6\}$ of $(\mathbb{Z}^6)^*$ (the dual) such that*

$$|y_{ijk} - x(u_i, u_j, u_k)| < \epsilon$$

*for all $i < j < k \leq 5$.*

We proved similar statements for cases (3), (6) in [25] (the irrationality condition is slightly easier).

We list $X_F$'s for Problem 9(1) except for the cases (1), (29).

(2) $\mathbb{P}(V)$
(3) $\mathbb{P}(V)$
(5) $\mathrm{Gr}(3,6), \mathrm{Gr}(3,6), (\mathbb{Z}/2\mathbb{Z}) \setminus (\mathrm{Gr}(3,6) \times \mathrm{Gr}(3,6))$
(6) $\mathbb{P}(\mathrm{Sym}^2(\mathbb{Q}^7)^*)$
(7) $\mathbb{P}(\mathrm{Sym}^2(\mathbb{Q}^8)^*)$
(10) $\mathbb{P}(\mathrm{Sym}^2(\mathbb{Q}^3)^*), \mathbb{P}(\mathrm{Sym}^2(\mathbb{Q}^5)^*)$
(13) $(2m \neq n)$ $\mathrm{Gr}(2m, 2n), \mathbb{P}(\wedge^2(\mathbb{Q}^{2m})^*)$
(13) $(2m = n)$ $(\mathbb{Z}/2\mathbb{Z}) \setminus (\mathrm{Gr}(2m, 4m) \times \mathrm{Gr}(2m, 4m)), \mathbb{P}(\wedge^2(\mathbb{Q}^{2m})^*)$
(14) the same as (5)
(15) $(2m \neq n)$ $\mathrm{Gr}(m, n), \mathbb{P}(\wedge^2(\mathbb{Q}^m)^*)$
(15) $(2m = n)$ $(\mathbb{Z}/2\mathbb{Z}) \setminus (\mathrm{Gr}(m, 2m) \times \mathrm{Gr}(m, 2m)), \mathbb{P}(\wedge^2(\mathbb{Q}^m)^*)$
(16) $\mathbb{P}(V)$
(18) $\mathrm{Gr}(3, 7), \mathbb{P}(\mathrm{Sym}^2(\mathbb{Q}^3)^*)$
(19) $\mathbb{P}^8$
(20) $\mathrm{Gr}(3, 10)$
(21) $(\mathbb{Z}/2\mathbb{Z}) \setminus (\mathrm{Gr}(5, 10) \times \mathrm{Gr}(5, 10)), \mathrm{Gr}(4, 10), \mathbb{P}^9, \mathbb{P}(\mathrm{Sym}^2(\mathbb{Q}^3)^*)$
(22) $\mathrm{Gr}(6, 12), \mathrm{Gr}(6, 12), (\mathbb{Z}/2\mathbb{Z}) \setminus (\mathrm{Gr}(6, 12) \times \mathrm{Gr}(6, 12))$
(23) the same as (22)
(24) $(\mathbb{Z}/2\mathbb{Z}) \setminus (\mathrm{Gr}(7, 14) \times \mathrm{Gr}(7, 14))$
(25) $\mathbb{P}(V)$
(27) $\mathbb{P}(V)$
(29) ?

We probably need some explanation. For (5), we can associate two 3-dimensional subspaces of a 6-dimensional space (which is the standard representation of $\mathrm{GL}(6)$) to any point in $V^{\mathrm{ss}}$. For this, the reader should see [25]. For (10), we can associate quadratic forms on two vector spaces $V_1, V_2$ of dimensions three and five (the prehomogeneous vector space (5) is considered as $\wedge^2 V_2 \otimes V_1$) to any point in $V^{\mathrm{ss}}$. For (18), we can associate a 3-dimensional subspace of a 7-dimensional space (which is the vector representation of $\mathrm{GSpin}(7)$) to any point in $V^{\mathrm{ss}}$. For (19), $\mathbb{P}^8$ is the projective space associated with the vector representation of $\mathrm{GSpin}(9)$. For (21), we can associate a 4-dimensional subspace, a 5-dimensional subspace of a 10-dimensional space (which is the vector representation of $\mathrm{GSpin}(10)$), and a quadratic form in three variables to any point in $V^{ss}$. By considering combinations of these subspaces, we get varieties as above. For (20), (22), (24), the numbers 10, 12, 14 are the dimensions of the vector representations. For (13) $(2m = n)$, (15) $(2m = n)$, (20), (24), we do not have to consider $\mathrm{Gr}(2m, 4m), \mathrm{Gr}(m, 2m), \mathrm{Gr}(5, 10), \mathrm{Gr}(7, 14)$ as we do for (5), (14), (22).

Consider Problem 9(2). Let $(G, V)$ be a prehomogeneous vector space in (1)–(29), and $\Delta(x)$ the relative invariant polynomial of the minimum degree. Suppose we carry out Problem 9(2). Then we conjecture the following application.

**Conjecture 12** If $g \in \mathrm{GL}(V)_\mathbb{R}$ is sufficiently irrational, the values of $\Delta(g^{-1}x)$ on the set of primitive integer points $x$ are dense in $\mathbb{R}$.

Note that if $f(x)$ is a polynomial on $V$, the natural action of $g$ on $f$ is given by $f(g^{-1}x)$. Let $d$ be the degree of $\Delta(x)$. Then since $(G, V)$ is a prehomogeneous vector space, $\Delta(x)$ is the unique $H_\mathbb{C}$ invariant of degree $d$ up to a constant. Therefore, we may choose $X_F = \mathbb{P}(\mathrm{Sym}^d V^*)$ and $w_F = [\Delta] \in \mathbb{P}(\mathrm{Sym}^d V^*)$ (the class of $\Delta$) for $F = H_{1\mathbb{C}}$ for Problem 9(2).

We will carry out Problem 9(2) and Conjecture 12 for the case (4) in the classification here. In this case, $H = \mathrm{SL}(2)$ and $H_2 = \mathrm{SL}(4)$. If $g$ is in the kernel of the homomorphism $H \to H_2$, $g$ is a scalar matrix $tI_2$ and $t^3 = 1$. Since this is an element of $\mathrm{SL}(2)$, $t^2 = 1$ and this implies $t = 1$. So we may regard $H$ as a subgroup of $H_2$. Therefore, $H_1 = H = \mathrm{SL}(2)$ in this case.

We regard $V$ as the space of cubic forms in two variables $v = (v_1, v_2)$. For $x = x_0 v_1^3 + x_1 v_2^2 v_2 + x_2 v_1 v_2^2 + x_3 v_2^3$, $y = y_0 v_1^3 + y_1 v_2^2 v_2 + y_2 v_1 v_2^2 + y_3 v_2^3$, we define

$$B(x, y) = x_0 y_3 - \frac{1}{3} x_1 y_2 + \frac{1}{3} x_2 y_1 - x_3 y_0.$$

Then it is well known that the alternating bilinear form $B$ is invariant under the action of $H_1$. For example, this $B$ is used in [24, p. 154]. We denote the symplectic group with respect to $B$ by $\mathrm{Sp}(4)$.

**Lemma 13** If $H_{1\mathbb{C}} \subset F \subset H_{2\mathbb{C}}$ is a closed connected complex Lie subgroup, $F = \mathrm{SL}(2)_\mathbb{C}, \mathrm{Sp}(4)_\mathbb{C},$ or $\mathrm{SL}(4)_\mathbb{C}$.

*Proof.* Let $\mathfrak{h}_1, \mathfrak{f}, \mathfrak{h}_2$ be the Lie algebras of $H_{1\mathbb{C}}, F, H_{2\mathbb{C}}$. These are vector spaces over $\mathbb{C}$. Let $\Lambda$ be the fundamental dominant weight of the Lie algebra $\mathfrak{h}_1$. We use the notation $d\Lambda$ for the irreducible representation of $\mathfrak{h}_1$ with highest weight $d\Lambda$.

By considering weights, $\mathfrak{h}_2 = \mathfrak{h}_1 \oplus 4\Lambda \oplus 6\Lambda$ and $\mathfrak{h}_1 \cong 2\Lambda$. Let $U_1 = 4\Lambda, U_2 = 6\Lambda$. Let $\mathfrak{h}_3 = \mathrm{sp}(4)_\mathbb{C}$ be the Lie algebra of $\mathrm{Sp}(4)_\mathbb{C}$. Counting the dimension, $\mathfrak{h}_3 = \mathfrak{h}_1 \oplus U_2$. Since $\mathfrak{h}_3$ is a Lie subalgebra of $\mathfrak{h}_2$, $\mathfrak{h}_2$ is a representation of $\mathfrak{h}_3$. Since $\mathfrak{h}_3$ is semisimple, there exists a representation $U_3$ of $\mathfrak{h}_3$ such that $\mathfrak{h}_2 = \mathfrak{h}_3 \oplus U_3$. Now $U_3$ can be regarded as a representation of $\mathfrak{h}_1$. Since $\mathfrak{h}_3$ is isomorphic to $2\Lambda \oplus 6\Lambda$, $U_3$ must be the factor $4\Lambda$. So $U_1 = 4\Lambda$ is a representation of $\mathfrak{h}_3$. This implies that $[U_1, U_1]$ is a representation of $\mathfrak{h}_3$ also. Since $\mathfrak{h}_2$ is a simple group, $[U_1, U_1] \cap \mathfrak{h}_3 \neq 0$. Since $\mathfrak{h}_3$ is the adjoint representation of a simple group, it is irreducible, and therefore, $[U_1, U_1] \supset \mathfrak{h}_3$. So we can conclude that if $\mathfrak{f} \supset U_1$, $\mathfrak{f} = \mathfrak{h}_2$. This proves the lemma. □

The groups $\mathrm{SL}(2), \mathrm{SL}(4), \mathrm{Sp}(4)$ are all defined over $\mathbb{Q}$. It is known [24, p. 150] that the relative invariant polynomial of the lowest degree is

$$\Delta(x) = x_1^2 x_2^2 + 18 x_0 x_1 x_2 x_3 - 4 x_0 x_2^3 - 4 x_1^3 x_3 - 27 x_0^2 x_3^2$$

for $x = x_0 v_1^3 + x_1 v_1^2 v_2 + x_2 v_1 v_2^2 + x_3 v_2^3$. Therefore, for $F = \mathrm{SL}(2)_{\mathbb{C}}$, we may choose $X_F = \mathbb{P}(\mathrm{Sym}^4 V^*)$ and $w_F = [\Delta]$ (the class of $\Delta$). It is well known that the alternating bilinear form fixed by a simplectic group is a constant multiple of the corresponding alternating bilinear form. Therefore, for $F = \mathrm{Sp}(4)_{\mathbb{C}}$, we may choose $X_F = \mathbb{P}(\wedge^2 V^*)$ and $w_F = [B]$ (the class of $B$). This is the answer to Problem 9(2) in this case.

**Theorem 14** *Suppose $g \in \mathrm{GL}(4)_{\mathbb{R}}$, $g[\Delta] \notin \mathbb{P}(\mathrm{Sym}^4 V^*)_{\mathbb{Q}}$, and $g[B] \notin \mathbb{P}(\wedge^2 V^*)_{\mathbb{Q}}$. Then the values of $\Delta(g^{-1}x)$ on the set of primitive integer points $x$ are dense in $\mathbb{R}$.*

*Proof.* If $\det g < 0$, we choose $g' \in \mathrm{GL}(4)_{\mathbb{Z}}$ such that $\det g' = -1$. Then $g[\Delta] \notin \mathbb{P}(\mathrm{Sym}^4 V^*)_{\mathbb{Q}}$, $g[B] \notin \mathbb{P}(\wedge^2 V^*)_{\mathbb{Q}}$ if and only if $g'g[\Delta] \notin \mathbb{P}(\mathrm{Sym}^4 V^*)_{\mathbb{Q}}$, $g'g[B] \notin \mathbb{P}(\wedge^2 V^*)_{\mathbb{Q}}$ and the values of $\Delta(g^{-1}x)$ on the set of primitive integer points $x$ are dense in $\mathbb{R}$ if and only if the values of $\Delta(g^{-1}g'^{-1}x)$ on the set of primitive integer points $x$ are dense in $\mathbb{R}$. Therefore, we may assume that $\det g > 0$. Then we can write $g = tI_4 g''$ where $t \in \mathbb{R}^{\times}$ and $g'' \in \mathrm{SL}(4)_{\mathbb{R}}$. Since $g[\Delta] = g''[\Delta]$ and $g[B] = g''[B]$, we may assume $g \in \mathrm{SL}(4)_{\mathbb{R}}$. Note that the center of $\mathrm{SL}(4)$ acts on $\mathbb{P}(\mathrm{Sym}^4 V^*), \mathbb{P}(\wedge^2 V^*)$ trivially.

Let $\Gamma = \mathrm{SL}(4)_{\mathbb{Z}}$. Clearly, $H_{1\mathbb{R}}, H_{2\mathbb{R}}$ are connected in the classical topology. So, by Theorem 7, $(gH_{1\mathbb{R}}g^{-1})\Gamma$ is dense in $H_{2\mathbb{R}}$. This implies that $\Gamma(gH_{1\mathbb{R}}g^{-1})$ is dense in $H_{2\mathbb{R}}$ also. Let $r \in \mathbb{R}^{\times}$. Let $\mathrm{diag}(t_1, t_2, t_3, t_4)$ be the diagonal matrix with diagonal entries $t_1, t_2, t_3, t_4$. Then $h = \mathrm{diag}(-\frac{4}{r}, -\frac{r}{4}, 1, 1) \in \mathrm{SL}(4)_{\mathbb{R}}$ and

$$\Delta(h^{-1}(1,0,1,0)) = r.$$

There exist $h_1 \in gH_{1\mathbb{R}}g^{-1}$ and $h_2 \in \Gamma$ such that $h_2 h_1$ is close to $hg^{-1}$. Then

$$(h_2 h_1 g \Delta)(1,0,1,0) = (h_2 g \Delta)(1,0,1,0) = \Delta(g^{-1}h_2^{-1}(1,0,1,0))$$

is close to $h\Delta(1,0,1,0) = \Delta(h^{-1}(1,0,1,0)) = r$. Since $h_2^{-1}(1,0,1,0)$ is a primitive integer point, this proves the theorem. □

Theorem 14 gives us a 13-dimensional family of quartic forms in four variables.

Question 2 can be considered as Problem 9(1) for the prehomogeneous vector space (2) $n = 3$ or Problem 9(2) for the prehomogeneous vector space (2) $n = 2$. All the analytical difficulties are contained in Ratner's theorem and Problem 5 for this case is almost trivial (which of course does not mean that Question 2 is trivial).

As far as Problem 9(2) is concerned, besides the prehomogeneous vector space (4), we carried out this program for the representation $(\mathrm{SL}(2), \mathrm{Sym}^4 k^2)$ [27], and the adjoint representation of $\mathrm{SL}(3)$ [26].

Let $V$ be the representation space for both cases. For these cases, $\dim V = 5$ or 8. It is well known that there are a quadratic invariant $Q(x)$ and a cubic invariant $F(x)$ for both cases. Then we have the following theorem [26], [27].

**Theorem 15** *For the above two cases, if $g \in \mathrm{GL}(V)_{\mathbb{R}}$ and $Q(g^{-1}x)$ is irrational, the values of $F(g^{-1}x)$ at primitive integer points are dense in $\mathbb{R}$.*

These are a 22-dimensional and a 56-dimensional family of cubic forms in 5 and 8 variables respectively. For cubic forms, Pitman [13] proved a weaker form of the Oppenheim conjecture (in the sense that if $f(x)$ is the form and $\epsilon > 0$, there exists an integer point $x \neq 0$ such that $|f(x)| < \epsilon$). Still this weaker form requires $(1314)^{256}$ variables, and it is of interest to construct examples of forms with dense values at primitive integer points.

Consider Conjecture 12 for the prehomogeneous vector spaces (1)–(29). If the degree of $\Delta$ is two, the conjecture is not so interesting because it is contained in the Oppenheim conjecture. This applies to cases (16), (19), (25). If a prehomogeneous vector space $(G_1, V_1)$ is contained in another prehomogeneous vector space $(G_2, V_2)$ and the restriction of a relative invariant polynomial on $V_2$ is a relative invariant polynomial on $V_1$, Conjecture 12 for $(G_2, V_2)$ is not so interesting. This applies to cases (1), (5), (12), (17), (22), (23), (26), (27), (28), (29), because the cases (1), (27) reduce to the case (2), the cases (5), (22), (23), (29) reduce to the case (14), the cases (12), (28) reduce to the case (8), and the cases (17), (26) reduce to the case (15) $m = 2$. It is not entirely clear, but probably the case (3) reduces to the case (2), the case (9) reduces to the case (8), and the case (13) reduces to the case (15). This leaves us the following cases.

| case | (2) | (4) | (6) | (7) | (8) | (10) | (11) | (14) | (15) | (20) | (21) | (24) |
|---|---|---|---|---|---|---|---|---|---|---|---|---|
| $\deg \Delta$ | $n$ | 4 | 7 | 16 | 12 | 15 | 40 | 4 | $2m$ | 4 | 4 | 8 |

Of course it is more interesting if $\deg \Delta$ is relatively large compared to the number of variables. In particular, by considering the cases (8), (11), it may be possible to produce families of polynomials of degree 12 and 40 in 12 and 40 variables respectively with dense values on the set of primitive integer points.

Finally it is curious to see if one can construct families of polynomials with non-trivial moduli (in other words, they consist of infinitely many orbits) with dense values on the set of primitive integer points. For example, it is possible to do so using the representation $(\mathrm{SL}(2), \mathrm{Sym}^4 \mathbb{Q}^2)$, but we have yet to know if this is possible for the spaces of binary forms of higher degrees for example.

**Acknowledgement** The author would like to thank D. Witte for reading the manuscript carefully.

Department of Mathematics, Oklahoma State University, Stillwater OK 74078
email: yukie@math.okstate.edu